\documentclass{amsart}
\usepackage{amssymb}
\usepackage{graphicx}

\newtheorem{thm}{Theorem}
\newtheorem{cor}{Corollary}
\newtheorem{lem}{Lemma}

\newtheorem{rem}{Remark}

\newcommand{\A}{{\mathcal A}}

\newcommand{\U}{{\mathcal U}}
\newcommand{\es}{{\mathcal S}}

\newcommand{\IN}{{\mathbb N}}

\newcommand{\D}{{\mathbb D}}

\def\be{\begin{equation}}
\def\ee{\end{equation}}

\newcommand{\bee}{\begin{enumerate}}
\newcommand{\eee}{\end{enumerate}}

\newcommand{\blem}{\begin{lem}}
\newcommand{\elem}{\end{lem}}
\newcommand{\bthm}{\begin{thm}}
\newcommand{\ethm}{\end{thm}}
\newcommand{\bcor}{\begin{cor}}
\newcommand{\ecor}{\end{cor}}
\newcommand{\beg}{\begin{example}}
\newcommand{\eeg}{\end{example}}
\newcommand{\begs}{\begin{examples}}
\newcommand{\eegs}{\end{examples}}
\newcommand{\bdefe}{\begin{defin}}
\newcommand{\edefe}{\end{defin}}
\newcommand{\bprob}{\begin{prob}}
\newcommand{\eprob}{\end{prob}}
\newcommand{\bei}{\begin{itemize}}
\newcommand{\eei}{\end{itemize}}

\newcommand{\bcon}{\begin{conj}}
\newcommand{\econ}{\end{conj}}
\newcommand{\bcons}{\begin{conjs}}
\newcommand{\econs}{\end{conjs}}
\newcommand{\bprop}{\begin{propo}}
\newcommand{\eprop}{\end{propo}}
\newcommand{\br}{\begin{rem}}
\newcommand{\er}{\end{rem}}
\newcommand{\brs}{\begin{rems}}
\newcommand{\ers}{\end{rems}}
\newcommand{\bo}{\begin{obser}}
\newcommand{\eo}{\end{obser}}
\newcommand{\bos}{\begin{obsers}}
\newcommand{\eos}{\end{obsers}}
\newcommand{\bpf}{\begin{pf}}
\newcommand{\epf}{\end{pf}}
\newcommand{\ba}{\begin{array}}
\newcommand{\ea}{\end{array}}
\newcommand{\beq}{\begin{eqnarray}}
\newcommand{\beqq}{\begin{eqnarray*}}
\newcommand{\eeq}{\end{eqnarray}}
\newcommand{\eeqq}{\end{eqnarray*}}

\begin{document}
\bibliographystyle{amsplain}

\title[Zalcman and generalized Zalcman conjecture for the class ${\U}$]{Zalcman and generalized Zalcman conjecture for the class $\boldsymbol{\U}$}

\author[M. Obradovi\'{c}]{Milutin Obradovi\'{c}}
\address{Department of Mathematics,
Faculty of Civil Engineering, University of Belgrade,
Bulevar Kralja Aleksandra 73, 11000, Belgrade, Serbia}
\email{obrad@grf.bg.ac.rs}

\author[N. Tuneski]{Nikola Tuneski}
\address{Department of Mathematics and Informatics, Faculty of Mechanical Engineering, Ss. Cyril and Methodius
University in Skopje, Karpo\v{s} II b.b., 1000 Skopje, Republic of North Macedonia.}
\email{nikola.tuneski@mf.edu.mk}





\begin{abstract}
Function $f(z)=z+\sum_{n=2}^{\infty} a_n z^n$, normalized, analytic and univalent in the unit disk $\D=\{z:|z|<1\}$, belongs to the class $\U$. if, and only if,
\[ \left| \left(\frac{z}{f(z)}\right)^2 f'(z) -1\right|<1 \quad\quad  (z\in\D). \]
In this paper we prove the Zalcman and the generalized Zalcman conjecture for the class $\U$ and some values of parameters in the conjectures.
\end{abstract}

\maketitle

\medskip

\section{Introduction and preliminaries}

\medskip

Let $\mathcal{A}$ be the class of functions $f$ which are analytic  in the open unit disc $\D=\{z:|z|<1\}$ of the form
\be\label{e1}
f(z)=z+a_2z^2+a_3z^3+\cdots,
\ee
and let $\mathcal{S}$ be the subclass of $\mathcal{A}$ consisting of functions that are univalent in $\D$.

\medskip

Important milestone in study of univalent functions was the proof of the famous Bieberbach conjecture $|a_n|\le n$ for $n\ge2$ by Lewis de Branges  in 1985 \cite{Bra85}.
This ended an era, but a great many other problems concerning the coefficients $a_n$ remain open. One such is the Zalcman conjecture,
\[ |a_n^2-a_{2n-1}|\le (n-1)^2 \quad\quad (n\in \IN, n\ge2), \]
posed in 1960 and proven in 2014 by Krushkal (\cite{krushkal}) for the whole class $\es$ by using complex geometry of the universal Teichm\"{u}ller space. In 1999, Ma (\cite{ma}) proposed a generalized Zalcman conjecture,
\[ |a_m a_n-a_{m+n-1}|\le (m-1)(n-1) \quad\quad (m,n\in \IN, m\ge2, n\ge2),\]
which is still an open problem, closed by Ma  for the class of starlike functions and for the class of univalent functions with real coefficients. Ravichandran and Verma in \cite{ravi} closed it for the classes of starlike and convex functions of given order and for the class of functions with bounded turning.

\medskip

In this paper we study the generalized Zalcman conjecture for the class
\[\U=\left\{f\in\A:\left| \left(\frac{z}{f(z)}\right)^2 f'(z) -1\right|<1, \, z\in\D \right\}.\]
Functions from this class are proven to be univalent but do not follow the traditional patterns of other univalent functions. For example, they are not starlike which makes them interesting since the class of starlike functions is very wide. So, class $\U$  attracts significant attention in the past decades and an overview of the most valuable results is given in Chapter 12 from \cite{DTV-book}.

\medskip

Here we will prove of the generalized Zalcman conjecture for the class $\U$ and for the cases $m=2$, $n=3$; $m=2$, $n=4$; and $m=n=3$.

\medskip

We also give direct proof and sharpness of the inequality
\[ |a_n^p-a_2^{p(n-1)}|\le 2^{p(n-1)}-n^p\]
over the class $\U$ for the cases $n=4$, $p=1$ and $n=5$, $p=1$. This inequality introduced by Krushkal and proven for the whole class of univalent functions in \cite{krushkal}.

\medskip

For the study and the proofs we will use the following useful property of functions in $\U$.

\begin{lem}[\cite{OB-NT-2018}] \label{lem-1}
For each function $f$ in $\U$, there exists function  $\omega_1$, analytic in the unit disk,  such that $|\omega_1(z)|\le|z|<1$ and $|\omega_1'(z)|\le1$ for all $z\in\D$, with
\be\label{u1}
\frac{z}{f(z)} = 1-a_2z-z\omega_1(z).
\ee
Additionally, for $\omega_1(z) = c_1z+c_2z^2+\cdots$,
\[|c_1|\le 1,\quad\quad   |c_2|\le \frac12(1-|c_1|^2) \quad \mbox{and}\quad |c_3|\le \frac13 \left( 1-|c_1|^2-\frac{4|c_2|^2}{1+|c_1|} \right). \]
\end{lem}

\medskip

Let note that for functions $f$ from $\U$, of form \eqref{e1}, from Lemma \ref{lem-1} we have
  \[ z = [1-a_2z-z\omega_1(z)]\cdot f(z), \]
and after equating the coefficients,
\[
\begin{split}
a_3 &= c_1+a_2^2,\\
a_4 &= c_2+2a_2c_1+a_2^3,\\
a_5 &= c_3+2a_2c_2+c_{1}^{2}+3a_2^2 c_{1}+a_{2}^{4}.
\end{split}
\]

\medskip

\section{Zalcman conjecture for the class $\U$}

\medskip

We now give direct proof of the Zalcman conjecture for the class $\U$ for the cases when  $n=2$ and $n=3$.

\medskip
	
\begin{thm}\label{th-1}
Let $f\in\U$ be of the form \eqref{e1}. Then
\begin{itemize}
  \item[$(i)$] $|a_2^2-a_3|\le1$;
  \item[$(ii)$] $|a_3^2-a_5|\le4$.
\end{itemize}
These inequalities are sharp with equality for the Koebe function $k(z)=\frac{z}{(1-z)^2}=z+\sum_{n=2}nz^n$ and its rotations.
\end{thm}

\medskip

\begin{proof}$ $
\begin{itemize}
  \item[$(i)$] From $a_3 = c_1+a_2^2$ we have $|a_2^2-a_3| = |-c_1|\le1$.
    \item[$(ii)$] From the Bieberbach conjecture, $|a_3| = |c_1+a_2^2|\le 3$, and further calculations show that
  \[
\begin{split}
\left|a_3^2-a_5 \right|
&= \left|c_3+2a_2c_2+a_2^2c_1 \right| \\
&= \left|c_3+2a_2c_2 - c_1^2+c_1(c_1 +a_2^2) \right|\\
&\le |c_3|+2|a_2||c_2| + |c_1|^2+|c_1||c_1 +a_2^2| \\
&\le |c_3|+2|a_2||c_2| + |c_1|^2+ 3|c_1| \\
&\le \frac13 \left( 1-|c_1|^2-\frac{4|c_2|^2}{1+|c_1|} \right) + 4|c_2| + |c_1|^2+ 3|c_1|\\
&:= f_1(|c_1|,|c_2|),
\end{split}
\]
where
\[ f_1(x,y) = \frac13 \left( 1-x^2-\frac{4y^2}{1+x} \right) + 4y + x^2+ 3x, \]
$0\le x=|c_1|\le1$ and $0\le y=|c_2|\le \frac12(1-x^2)$, i.e., $(x,y)\in G:=[0,1]\times[0,(1-x^2)/2]$.

\medskip

Since, $\frac{\partial f_1}{\partial y}(x,y)  = \frac43 \left( \frac{y}{1+x} \right)^2+\frac43x+3 >0$ for all $(x,y)\in G$, we have that there are no singular points in the interior of $G$ and $f_1$ attains its maximum on the boundary of $G$.

\medskip

Further, for $x=0$ we have $0\le y\le \frac12$ and $f_1(0,y)=\frac13(1-4y^2)+4y\le2$.

\medskip
Also, for $0\le x\le1$ and $y=0$ we have $f_1(x,0)=\frac13(1-x^2)+x^2+3x\le4$.

\medskip
Finally, for $0\le x\le1$ and $y=\frac12(1-x^2)$ we have $f_1(x,\frac12(1-x^2))=2+\frac{10}{3}x-x^2-\frac13x^3 \le4$, since the last function is an increasing one on $[0,1]$.
\end{itemize}
\end{proof}  

\medskip

\section{Generalized Zalcman conjecture for the class $\U$}

\medskip

In this section we give direct proof of the generalized Zalcman conjecture for the class $\U$ for the cases $m=2$, $n=3$; and $m=2$, $n=4$.

\medskip

\begin{thm}
Let $f\in\U$ be of the form \eqref{e1}. Then
\begin{itemize}
  \item[$(i)$] $|a_2a_3-a_4|\le2$;
  \item[$(ii)$] $|a_2a_4-a_5|\le3$.
\end{itemize}
These inequalities are sharp with equality for the Koebe function $k(z)=\frac{z}{(1-z)^2}=z+\sum_{n=2}nz^n$ and its rotations.
\end{thm}

\medskip

\begin{proof}$ $
\begin{itemize}
  \item[$(i)$] In this case we have
  \[
\begin{split}
|a_2a_3-a_4|
&= |c_2+a_2c_1| \le  |c_2|+|a_2||c_1| \le |c_2|+2|c_1| \\
&\le \frac12 (1-|c_1|^2)+ 2|c_1| \le \frac12 (1-|c_1|^2+4|c_1| )\le 2.
\end{split}
\]
  \item[$(ii)$] In a similar way as in the proof of Theorem \ref{th-1}($ii$),
    \[
\begin{split}
\left|a_4a_2-a_5 \right|
&= \left|c_3 + a_2c_2 + a_2^2c_1 + c_1^2 \right| \\
&\le |c_3| + |a_2||c_2| + |c_1| |a_2^2+c_1| \\
&\le |c_3| + |a_2||c_2| + 3|c_1| \\
&\le \frac13 \left( 1-|c_1|^2-\frac{4|c_2|^2}{1+|c_1|} \right) + 2|c_2| + 3|c_1|\\
&:= f_2(|c_1|,|c_2|),
\end{split}
\]
where
\[ f_2(x,y) = \frac13 \left( 1-x^2-\frac{4y^2}{1+x} \right) + 2y + 3x, \]
$0\le x=|c_1|\le1$ and $0\le y=|c_2|\le \frac12(1-x^2)$, i.e., $(x,y)\in G:=[0,1]\times[0,(1-x^2)/2]$.
\medskip

Again, $\frac{\partial f_2}{\partial y}(x,y)  = \frac43 \left( \frac{y}{1+x} \right)^2-\frac23x+3 >0$ for all $(x,y)\in G$, so $f_2$ attains its maximum on the boundary of $G$.

\medskip

The conclusion follows since on the edges of $G$ we have:
\begin{itemize}
  \item[-] $x=0$, $0\le y\le \frac12$ and $f_2(0,y)=\frac13(1-4y^2)+2y\le1$;
  \item[-] $y=0$,  $0\le x\le1$ and $f_2(x,0)=\frac13(1-x^2)+3x\le3$;
  \item[-] $y=\frac12(1-x^2)$, $0\le x\le1$ and $f_2(x,\frac12(1-x^2))=1+\frac{10}{3}x-x^2-\frac13x^3 \le3$.
\end{itemize}
\end{itemize}
\end{proof}

\medskip

\section{Krushkal inequality for the class $\U$}

\medskip

In this section we give direct proof of the Krushkal inequality for the class $\U$ in the cases when $n=4$, $p=1$ and $n=5$, $p=1$.

\medskip

\begin{thm}
Let $f\in\U$ be of the form \eqref{e1}. Then
\begin{itemize}
  \item[$(i)$] $|a_4-a_2^3|\le4$;
  \item[$(ii)$] $|a_5-a_2^4|\le11$.
\end{itemize}
These inequalities are sharp with equality for the Koebe function $k(z)=\frac{z}{(1-z)^2}=z+\sum_{n=2}nz^n$ and its rotations.
\end{thm}

\medskip

\begin{proof}$ $
\begin{itemize}
  \item[$(i)$] It is easy to verify that
   \[
\begin{split}
|a_4-a_2^3| &= |c_2+2a_2c_1| \le  |c_2|+2|a_2||c_1| \\
&\le \frac12 (1-|c_1|^2)+ 4|c_1| = \frac12\left(1+8|c_1|-|c_1|^2\right)\le 4.
\end{split}
\]
\medskip
  \item[$(ii)$] We will again use that $|a_3| = |c_1+a_2^2|\le 3$ and receive
   \[
\begin{split}
|a_5-a_2^4| &= |c_3+2a_2c_2+c_1^2+3a_2^2c_1| \\
& = |c_3+2a_2c_2-2c_1^2+3c_1(c_1+a_2^2)| \\
 &\le  |c_3|+2|a_2||c_2|+2|c_1|^2+9|c_1|  \\
&\le \frac13 \left( 1-|c_1|^2-\frac{4|c_2|^2}{1+|c_1|} \right) + 4|c_2| + 2|c_1|^2 + 9|c_1|\\
&:= g(|c_1|,|c_2|),
\end{split}
\]
where
\[ g(x,y) = \frac13 \left( 1-x^2-\frac{4y^2}{1+x} \right) + 4y + 2x^2+9x, \]
$0\le x=|c_1|\le1$ and $0\le y=|c_2|\le \frac12(1-x^2)$, i.e., $(x,y)\in G:=[0,1]\times[0,(1-x^2)/2]$.
\medskip

Since, $\frac{\partial g}{\partial y}(x,y)  = \frac{10}{3}x + \frac43 \left( \frac{y}{1+x} \right)^2 + 9 >0$ for all $(x,y)\in G$, so $g$ has no critical points in the interior of $G$ and attains its maximum on the boundary:
\begin{itemize}
  \item[-] $x=0$, $0\le y\le \frac12$ and $g(0,y)=\frac13(1+12y-4y^2)\le2$;
  \item[-] $y=0$,  $0\le x\le1$ and $g(x,0)=\frac53x^2+9x+\frac13\le11$;
  \item[-] $y=\frac12(1-x^2)$, $0\le x\le1$ and $g(x,\frac12(1-x^2))=2+\frac{28}{3}x-\frac13x^3 \le11$.
\end{itemize}
\medskip

The statement ($ii$) follows directly.
\end{itemize}
\end{proof}


\medskip


\begin{thebibliography}{9}

\bibitem{Bra85} L. De Branges, \emph{A proof of the Bieberbach conjecture}, Acta Math. {\bf 154} (1985), no. 1--2, 137--152.



%
%
%
%

\bibitem{krushkal}
S.L. Krushkal, A short geometric proof of the Zalcman and Bieberbach conjectures, preprint, 	arXiv:1408.1948.

\bibitem{ma}
 W. Ma, Generalized Zalcman conjecture for starlike and typically real functions, J. Math. Anal. Appl., 234 (1999), pp. 328-–339




\bibitem{OB-NT-2018}
M. Obradovi\'{c}, N. Tuneski, Some properties of the class $\U$, Annales. Universitatis Mariae Curie-Skłodowska. Sectio A - Mathematica, Vol. 73 No.1, (2019), 49--56.


\bibitem{ravi}
V. Ravichandran, S. Verma,  Generalized Zalcman conjecture for some classes of analytic functions. J. Math. Anal. Appl. 450 (2017), no. 1, 592--605.

\bibitem{DTV-book}
D.K. Thomas, N. Tuneski, A. Vasudevarao, Univalent Functions: A Primer, \emph{De Gruyter Studies in Mathematics} {\bf 69}, De Gruyter, Berlin, Boston, 2018.



\end{thebibliography}
\end{document}